\documentclass{mrlart3}
\usepackage{amsmath,amssymb,amsfonts,amscd,amsthm}
\usepackage[dvips]{graphicx}
\usepackage{psfrag,rotating}
\usepackage[colorlinks=true,%
		linkcolor=blue,%
		citecolor=blue]{hyperref}

\newcommand{\cL}{\mathcal{L}}

\newcommand{\hZ}{\hat{Z}}

\newcommand{\fS}{\mathfrak{S}}
\newcommand{\fP}{\mathfrak{P}}
\newcommand{\Tr}{\mathrm{Tr}}

  \def\volno{00}
  \def\yrno{0000}
  \def\issueno{00}
  \def\lpageno{0000} 
\title{On a proof of the Labastida-Mari\~no-Ooguri-Vafa Conjecture}

\author{Kefeng Liu}
  \address{Center of Mathematical Science \\
    Zhejiang University \\
    Box 310027 \\
    Hangzhou, China}
  \email{liu@cms.zju.edu.cn}

  \address{Department of Mathematics \\
    UCLA \\
    Box 951555 \\
    Los Angeles, CA, 90095-1555}
  \email{liu@math.ucla.edu}

\author{Pan Peng}
  \address{Department Of Mathematics \\
    University of Arizona \\
    617 N. Santa Rita Ave \\
    Tucson, AZ, 85721.}
  \email{ppeng@math.arizona.edu}

\date{}

\numberwithin{equation}{section} \theoremstyle{plain}
\newtheorem{Theorem}{Theorem}
\newtheorem{theorem}{Theorem}[section]

\newtheorem{corollary}[theorem]{Corollary}
\newtheorem{definition}[theorem]{Definition}
\newtheorem{lemma}[theorem]{Lemma}
\newtheorem{proposition}[theorem]{Proposition}
\newtheorem{remark}[theorem]{Remark}

\DeclareMathOperator{\Aut}{Aut} \DeclareMathOperator{\Ord}{Ord}
 \DeclareMathOperator{\End}{End}
\DeclareMathOperator{\tr}{tr} \DeclareMathOperator{\id}{id}

\newtheorem*{lmovconj}{Conjecture (LMOV)}

\DeclareMathOperator{\lk}{lk}
\newcommand{\chR}{\check{\mathcal{R}}}
\newcommand{\fsl}{\mathfrak{sl}}

\newdimen\tableauside\tableauside=1.0ex
\newdimen\tableaurule\tableaurule=0.4pt
\newdimen\tableaustep
\def\phantomhrule#1{\hbox{\vbox to0pt{\hrule height\tableaurule width#1\vss}}}
\def\phantomvrule#1{\vbox{\hbox to0pt{\vrule width\tableaurule height#1\hss}}}
\def\sqr{\vbox{%
  \phantomhrule\tableaustep
  \hbox{\phantomvrule\tableaustep\kern\tableaustep\phantomvrule\tableaustep}%
  \hbox{\vbox{\phantomhrule\tableauside}\kern-\tableaurule}}}
\def\squares#1{\hbox{\count0=#1\noindent\loop\sqr
  \advance\count0 by-1 \ifnum\count0>0\repeat}}
\def\partition#1{\vcenter{\offinterlineskip
  \tableaustep=\tableauside\advance\tableaustep by-\tableaurule
  \kern\normallineskip\hbox
    {\kern\normallineskip\vbox
      {\gettableau#1 0 }%
     \kern\normallineskip\kern\tableaurule}%
  \kern\normallineskip\kern\tableaurule}}
\def\gettableau#1 {\ifnum#1=0\let\next=\null\else
  \squares{#1}\let\next=\gettableau\fi\next}

\begin{document}
\maketitle
\begin{abstract}
 We outline a proof of a remarkable conjecture of Labastida-Mari\~no-Ooguri-Vafa
about certain new algebraic structures of quantum link invariants
and the integrality of infinite family of new topological
invariants. Our method is based on the cut-and-join analysis and a
special rational ring characterizing the structure of the
Chern-Simons partition function.
\end{abstract}

\section{Introduction}

For decades, we have witnessed the great development of string theory and its
powerful impact on the development of mathematics. There have been a lot of
marvelous results revealed by string theory, which deeply relate different
aspects of mathematics. All these mysterious relations are connected by a core
idea in string theory called ``duality". It was found that string theory on
Calabi-Yau manifolds provided new insight in geometry of these spaces. The
existence of a topological sector of string theory leads to a simplified model
in string theory, the topological string theory.

A major problem in topological string theory is how to compute Gromov-Witten
invariants. There are two major methods widely used: mirror symmetry in
physics and localization in mathematics. Both methods are effective
when genus is low while having trouble in dealing with higher genera due to
the rapidly growing complexity during computation. However, when the
target manifold is Calabi-Yau threefold, large $N$ Chern-Simons/topological
string duality opens a new gate to a complete solution of computing
Gromov-Witten invariants at all genera.

The study of large $N$ Chern-Simons/topological string duality was originated in
physics by an idea that gauge theory should have a string theory explanation. In
1992, Witten \cite{W2} related topological string theory of $T^\ast M$ of a
three dimensional manifold $M$ to Chern-Simons gauge theory on $M$. In 1998,
Gopakumar and Vafa \cite{GV} conjectured that, at large $N$, open topological
A-model of $N$ D-branes on $T^{\ast }S^{3}$ is dual to topological closed string
theory on resolved conifold $\mathcal{O}(-1)\oplus\mathcal{O}(-1)\rightarrow
\mathbb{P}^1$. Later, Ooguri and Vafa \cite{OV} showed a picture on how to
describe Chern-Simons invariants of a knot by open topological string theory on
resolved conifold paired with lagrangian associated with the knot.

Though large $N$ Chern-Simons/topological string duality still remains open,
there have been a lot of progress in this direction demonstrating the power of
this idea. Even for the simplest knot, the unknot, Mari\~{n}o-Vafa formular
\cite{MV, LLZ1} gives a beautiful closed formula for Hodge integral up to three
Chern classes of Hodge bundle. Furthermore, using topological vertex theory
\cite{AKMV, LLZ2, LLLZ}, one is able to compute Gromov-Witten invariants of any
toric Calabi-Yau threefold by reducing the computation to a gluing algorithm of
topological vertex. This thus leads to a closed formula of topological string
partition function, a generating function of Gromov-Witten invariants, in all
genera for any toric Calabi-Yau threefolds.

On the other hand, after Jones' famous work on polynomial knot
invariants, there had been a series of polynomial invariants
discovered (for example, \cite{J, HOMFLY, K}), the generalization of
which was provided by quantum group theory \cite{T} in mathematics
and by Chern-Simons path integral with the gauge group $SU(N)$
\cite{W1} in physics.

Based on the large $N$ Chern-Simons/topological string duality, Ooguri and Vafa
\cite{OV} reformulated knot invariants in terms of new integral invariants
capturing the spectrum of M2 branes ending on M5 branes embedded in the resolved
conifold. Later, Labastida, Mari\~{n}o and Vafa \cite{LMV,LM} refined the
analysis of \cite{OV} and conjectured the precise integrality structure for open
Gromov-Witten invariants. This conjecture predicts a remarkable new algebraic
structure for the generating series of general link invariants and the
integrality of infinite family of new topological invariants. In string theory,
this is a striking example that two important physical theories, topological
string theory and Chern-Simons theory, exactly agree up to all orders. In
mathematics this conjecture has interesting applications in understanding the
basic structure of link invariants and three manifold invariants, as well as the
integrality structure of open Gromov-Witten invariants. Recently, X.S. Lin and
H. Zheng \cite{LZ} verified LMOV conjecture in several lower degree cases for
some torus links.

In this paper, we describe an outline of the proof of
Labastida-Mari\~{n}o-Ooguri-Vafa conjecture for any link. The
details of the proofs are given in \cite{LP}.

\section{The Labastida-Mari\~no-Ooguri-Vafa conjecture}

\subsection{Quantum group invariants of links}

Let $\mathcal{L}$ be a link with $L$ components $\mathcal{K}_{\alpha
}$, $\alpha=1,\ldots,L$, represented by the closure of an element of
braid group $\mathcal{B}_{m}$. We associate to each
$\mathcal{K}_\alpha$ an irreducible representation $R_\alpha$ of
quantized universal enveloping algebra $U_q
(\fsl(N,\mathbb{C}))$\footnote{Later, we simply write
$U_q(\fsl_N)$.}, labeled by its highest weight $\Lambda_{\alpha }$.
Denote the corresponding module by $V_{\Lambda_\alpha}$. The $j$-th
strand in the braid will be associated with the irreducible module
$V_{j}=V_{\Lambda _{\alpha }}$, if this strand belongs to the
component $\mathcal{K}_{\alpha }$. The braiding is defined through
the following \emph{universal $R$-matrix} of $U_q(\fsl_N)$
\[
 \mathcal{R}
=q^{\frac{1}{2}\sum_{i,j}C_{ij}^{-1}H_i\otimes H_j}
 \prod_{\textrm{positive root }\beta} \exp_q
 [( 1-q^{-1}) E_\beta\otimes F_\beta]\,.
\]
Here $\{H_i,E_i,F_i\}$ are the generators of $U_q(\fsl_N)$,
 $( C_{ij})$ is the Cartan matrix and
\[
 \exp_{q}(x)
=\sum_{k=0}^{\infty }q^{\frac{1}{4}k(k+1)}\frac{x^{k}}{\{k\}_{q}!}\,,
\]
where
\begin{align*}
 \{k\}_q = \frac{q^{-k/2}-q^{k/2}}{q^{-1/2}-q^{1/2}}\,,
&&
 \{k\}_q!= \prod_{j=1}^k \{j\}_q\,.
\end{align*}
Define \emph{braiding} by $\check{\mathcal{R}}=P_{12}\mathcal{R}$, where
$P_{12}(v\otimes w)=w\otimes v$.

\begin{figure}[!ht]
\begin{align*}
    \psfrag{A}[][]{$U$}
    \psfrag{B}[][]{$V$}
    \psfrag{P}[][]{$\chR_{U,V}$}
    \includegraphics[width=1.5cm]{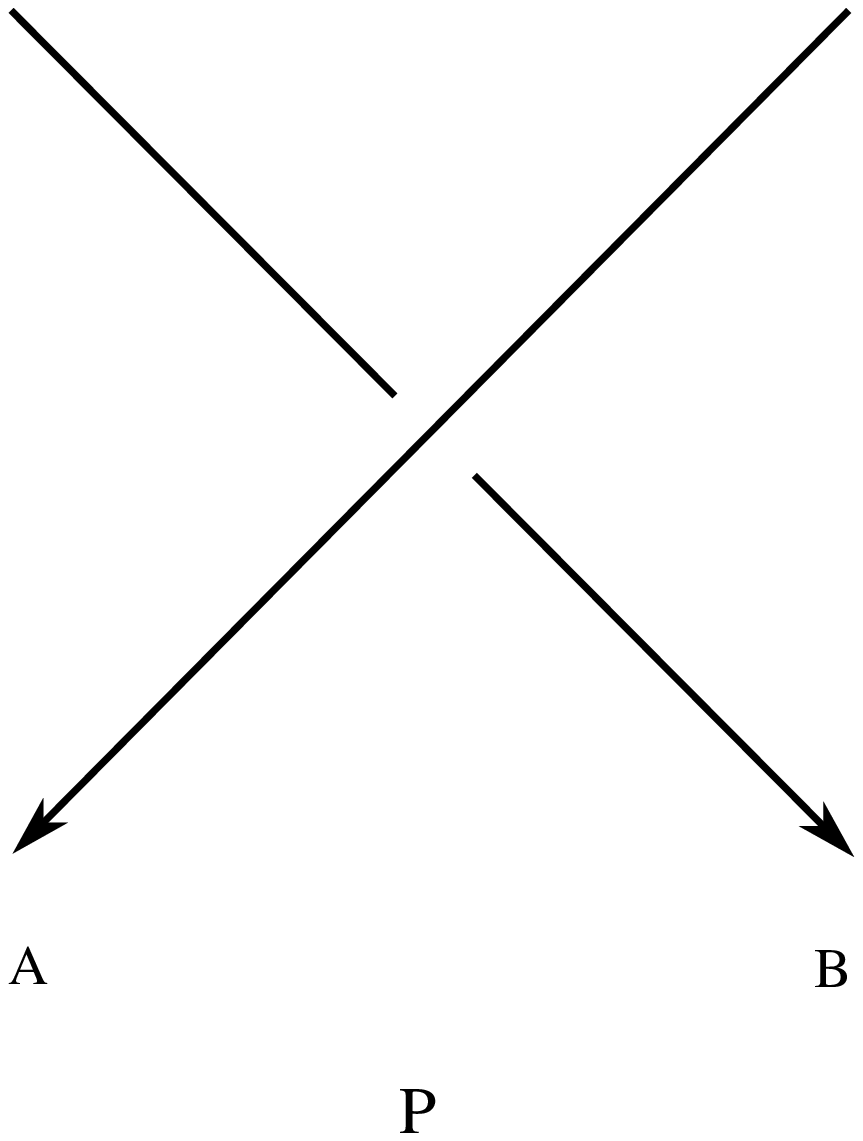}
  &&
    \psfrag{A}[][]{$U$}
    \psfrag{B}[][]{$V$}
    \psfrag{N}[][]{$\chR^{-1}_{V,U}$}
    \includegraphics[width=1.5cm]{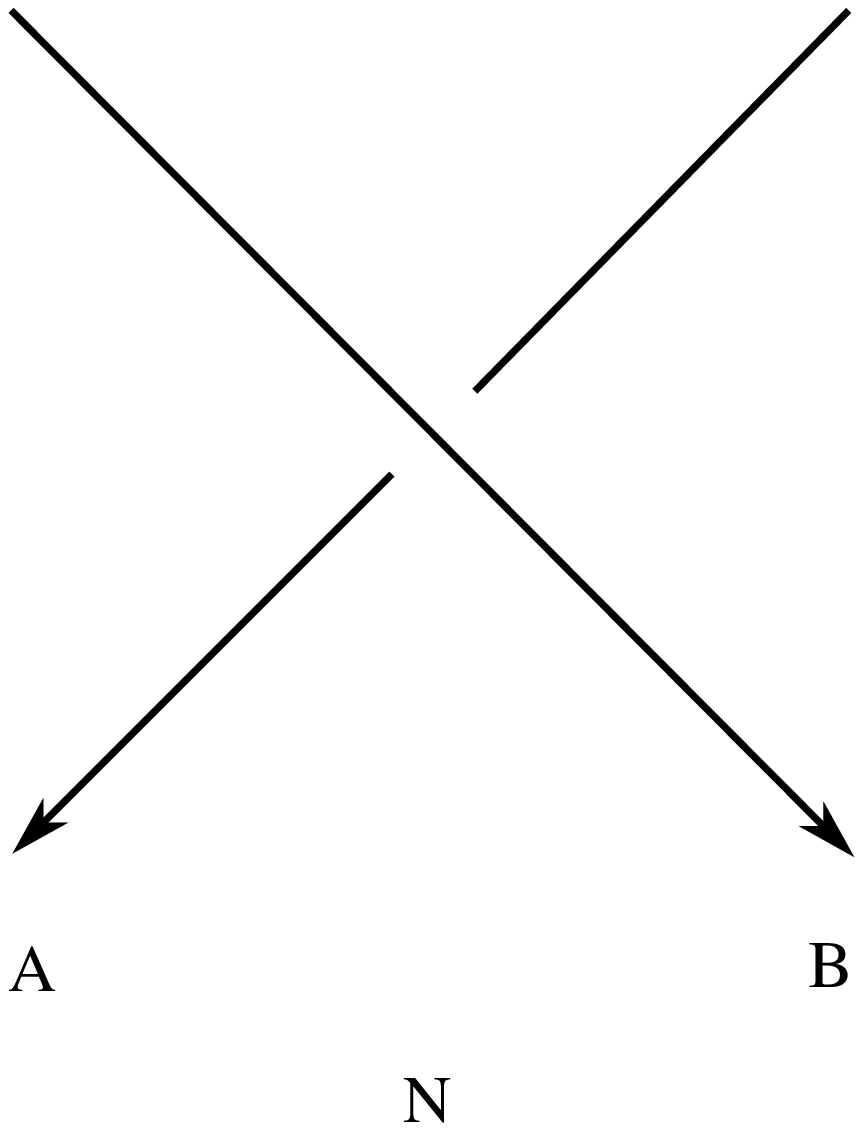}
\end{align*}
\caption{Assign crossing by $\chR$.}\label{fig: crossings}
\end{figure}

Now for a given link $\mathcal{L}$ of $L$ components, one chooses a closed braid
representative in braid group $\mathcal{B}_m$ whose closure is $\mathcal{L}$. In
the case of no confusion, we also use $\mathcal{L}$ to refer its braid
representative in $\mathcal{B}_m$. We will assign each crossing by the
braiding as follows. Let $U$, $V$ be two $U_q(\fsl_N)$-modules labeling two
outgoing strands of the crossing, the braiding $\check{R}_{U,V}$ (resp.
$\check{R}_{V,U}^{-1}$) is assigned as in Figure \ref{fig: crossings}.

The above assignment will give a representation of $\mathcal{B}_m$ on
$U_q(\mathfrak{g})$-module $V_1 \otimes \cdots \otimes V_m$. Namely, for any
generator,
$\sigma_i \in\mathcal{B}_m$, define\footnote{In the case of $\sigma_i^{-1}$, use
$\chR^{-1}_{V_{i+1}, V_i}$ instead.}
\begin{align*}
 h(\sigma_i)
= \id_{V_1} \otimes \cdots \otimes
 \chR_{V_i, V_{i+1}} \otimes \cdots \otimes id_{V_N}\,.
\end{align*}
Therefore, any link $\mathcal{L}$ will provide an isomorphism
\[
 h(\mathcal{L})
\in \End_{U_q (\fsl_N)}(V_{1}\otimes \cdots \otimes V_{m}) \, .
\]

Let $K_{2\rho}$ be the enhancement of $\check{\mathcal{R}}$ in the sense of
\cite{Rosso-Jones93}, where $\rho $ is the half-sum of all positive roots
of $\fsl_N$. The irreducible representation $R_{\alpha }$ is labeled by the
corresponding partition $A^\alpha$.

\begin{definition}\label{def: quantum group invariant}
Given $L$ labeling partitions $A^1,\ldots,A^L$, the quantum group invariant of
 $\mathcal{L}$ is defined as follows:
\[
    W_{(A^{1},...,A^{L})}(\mathcal{L})
=q^{d(\mathcal{L)}}\tr_{V_1\otimes\cdots\otimes V_{m}}
    (K_{2\rho}\circ h(\mathcal{L})) \, ,
\]
where
\[
 d(\mathcal{L})
=-\frac{1}{2}\sum_{\alpha =1}^{L}\omega (\mathcal{K}_\alpha)
  (\Lambda _{\alpha },\Lambda _{\alpha }+2\rho )
 +\frac{1}{N}\sum_{\alpha<\beta }
  \lk(\mathcal{K}_{\alpha },\mathcal{K}_{\beta })
  |A^{\alpha }|\cdot | A^{\beta }| \,,
\]
and $\lk(\mathcal{K}_\alpha,\mathcal{K}_\beta)$ is the linking number of
components $\mathcal{K}_\alpha$ and $\mathcal{K}_\beta$. A substitution of
$t=q^N$ is used to give a two-variable framing independent link invariant.
\end{definition}

\subsection{Labastida-Mari\~{n}o-Ooguri-Vafa conjecture}
\label{subsec: LMOV conjecture}

Let $\mathcal{L}$ be a link with $L$ components and $\mathcal{P}$ be the set of
all partitions. The Chern-Simons partition function of $\mathcal{L}$ is a
generating function of quantum group invariants of links given by
\begin{align}\label{eqn: definition of Chern-Simons partition function in sec0}
 Z_{\textrm{CS}}(\mathcal{L};\, q,t)
=\sum_{\vec{A}\in \mathcal{P}^L}W_{\vec{A}}(\mathcal{L};\, q,t)
    \prod_{\alpha=1}^L s_{A^\alpha}(x^\alpha)
\end{align}
for any arbitrarily chosen sequence of variables
\[
 x^\alpha =(x^\alpha_1,x^\alpha_2,\ldots,)\,.
\]
In \eqref{eqn: definition of Chern-Simons partition function in sec0},
$\vec{A}=(A^1,\ldots,A^L)\in \mathcal{P}^L$ and $s_{A^\alpha}(x^\alpha)$ is the
Schur function.

Free energy is defined to be
\[
 F=\log Z_{\textrm{CS}}\,.
\]
Use plethystic exponential, one can obtain
\begin{equation}\label{eqn: formula of f_A}
  F
= \sum_{d=1}^{\infty }\sum_{\vec{A}\neq 0}
    \frac{1}{d}f_{\vec{A}}(q^{d},t^{d})
    \prod_{\alpha =1}^{L}
    s_{A^{\alpha }}\big( ( x^\alpha)^{d} \big)\,,
\end{equation}
where
\[
 ( x^{\alpha }) ^{d}
= \big((x_{1}^\alpha)^d, (x_{2}^\alpha)^d,\ldots \big) \,.
\]

Based on the duality between Chern-Simons gauge theory and topological string
theory, Labastida, Mari\~{n}o, Ooguri, Vafa conjectured that $f_{\vec{A}}$ have
the following highly nontrivial structures.

For any $A$, $B\in \mathcal{P}$, define the following function
\begin{equation}\label{eqn: formula of M_AB}
 M_{AB}(q)
=\sum_{\mu} \frac{\chi_A(C_\mu) \chi_B(C_\mu)}{\mathfrak{z}_\mu}
 \prod_{j=1}^{\ell(\mu)} (q^{-\mu_j/2}-q^{\mu_j/2})\, .
\end{equation}

\begin{lmovconj}\label{conj: LMOV}
For any $\vec{A}\in \mathcal{P}^{L}$,
\begin{itemize}
\item[$(\mathrm{i})$.] there exist $P_{\vec{B}}(q,t)$ for
$\vec{B}\in\mathcal{P}^L$, such that
\begin{equation}\label{eqn: f_A of P_B}
 f_{\vec{A}}(q,t)
=\sum_{ |B^\alpha|=|A^\alpha|} P_{\vec{B}}(q,t)
 \prod_{\alpha =1}^L M_{A^\alpha B^\alpha }(q).
\end{equation}
Furthermore, $P_{\vec{B}}(q,t)$ has the following expansion
\begin{equation}\label{eqn: P_B}
 P_{\vec{B}}(q,t)
=\sum_{g= 0}^\infty \sum_{Q\in \mathbb{Z}/2} N_{\vec{B};\,g,Q}
 (q^{-1/2}-q^{1/2})^{2g-2}t^{Q}\,.
\end{equation}

\item[$(\mathrm{ii})$.] $N_{\vec{B};\,g,Q}$ are integers.
\end{itemize}
\end{lmovconj}

\subsection{Related notations}

For a given link $\mathcal{L}$ of $L$ components, we will fix the following
notations in this paper. Given
$\lambda\in\mathcal{P}$, $\vec{A}=(A^1,\ldots,A^L)$, $\vec{\mu} =
(\mu^1,\ldots, \mu^L) \in \mathcal{P}^L$.
Let $x=(x^1,...,x^L)$ where $x^\alpha$ is a set of variables
\[
 x^\alpha=(x^\alpha_1,x^\alpha_2,\cdots)\,.
\]
The following notations will be used
throughout the paper.
\begin{align*}
 &
 [n]_q = q^{-\frac{n}{2}}-q^{\frac{n}{2}}\,,
    &&
 [\lambda]_q = \prod_{j=1}^{\ell(\lambda)} [\lambda_j]_q\,,
    &&
 \mathfrak{z}_{\vec{\mu}}=\prod_{\alpha =1}^{L}\mathfrak{z}_{\mu^\alpha}\,,
\\
 &
 |\vec{A}| = ( |A^{1}|,...,|A^{L}|)\,,
    &&
 \parallel\vec{A}\parallel =\sum_{\alpha=1}^{L}|A^\alpha|\,,
    &&
 \ell(\vec{\mu}) = \sum_{\alpha=1}^{L}\ell (\mu ^\alpha) \,,
\\
 &
 \vec{A}^{t} =\big( ( A^{1})^t,\ldots,(A^{L})^{t}\big)\,,
    &&
 \chi _{\vec{A}}\left( \vec{\mu}\right)=\prod_{\alpha =1}^{L}
    \chi _{A^{\alpha }}( C_{\mu ^\alpha})\,,
    &&
 s_{\vec{A}}(x)=\prod_{\alpha=1}^{L}s_{A^\alpha}(x^{\alpha }) \,.
\end{align*}

One can define an order on $\mathcal{P}^L$ lexicographically.
Therefore, one can generalize the concept of partition from the set
of all non-negative integers to $\mathcal{P}^L$. We denote by
$\mathcal{P(P}^L)$ the set of all partitions on $\mathcal{P}^L$.
Given $\Lambda\in \mathcal{P(P}^L)$, the following quantity
\[
 \theta_\Lambda
=\frac{(-1)^{\ell(\Lambda)-1} (\ell(\Lambda)-1)!}{|\Aut\Lambda|}
\]
plays an important role in the relationship of topological
string partition function and free energy.

Rewrite free energy as
\begin{equation}\label{eqn: def of F_mu}
F=\log Z=\sum_{\vec{\mu}\neq 0}F_{\vec{\mu}}p_{\vec{\mu}}\left( x%
\right) .
\end{equation}
Here in the similar usage of notation,
\[
 p_{\vec{\mu}}(x)=\prod_{\alpha=1}^L p_{\mu^\alpha}(x^\alpha)\,.
\]
We also rewrite Chern-Simons partition function as
\[
 Z_{\mathrm{CS}}(\mathcal{L})
=1+\sum_{\vec{\mu}\neq 0}Z_{\vec{\mu}}p_{\vec{\mu}}(x)\,.
\]

\section{Sketch of the proof}

In this section we will give an outline of the proof of
Labastida-Mari\~no-Ooguri-Vafa conjecture\footnote{Briefly, we call it LMOV
 conjecture.} which contains two parts in correspondence of (i) and (ii) in
the
description of LMOV conjecture: the \textbf{existence and integrality}.

\subsection{Existence of the algebraic structure}

The existence of the algebraic structure \eqref{eqn: P_B} includes
the following two steps:
\begin{itemize}
 \item
The symmetry of $q$ and $q^{-1}$ in $P_{\vec{B}}(q,t)$.

 \item
The pole structure of $P_{\vec{B}}$.
\end{itemize}

The symmetry of $q$ and $q^{-1}$ in $P_{\vec{B}}(q,t)$ can
be obtained from the lemma.

\begin{lemma}\label{lemma: symmetry}
$W_{\vec{A}^{t}}(q,t)=(-1)^{\parallel\vec{A}\parallel}W_{\vec{A}}(q^{-1},t)$.
\end{lemma}

To prove the existence of the pole structure, we will consider the following
framed generating series. Substitute
\[
 W_{\vec{A}}(\mathcal{L};\,q,t,\tau )
=W_{\vec{A}}(\mathcal{L};\,q,t)\cdot
    q^{\sum_{\alpha=1}^{L}\kappa_{A^\alpha}\tau/2}
\]
in the Chern-Simons partition function, we have the following framed partition
function
\[
 Z( \mathcal{L};\,q,t,\tau)
=1+\sum_{\vec{A}\neq 0}W_{\vec{A}}(\mathcal{L};\,q,t,\tau )\cdot
    s_{\vec{A}}(x) \,.
\]
Similarly, framed free energy
\[
 F(\mathcal{L};\,q,t,\tau)
= \log Z( \mathcal{L};\,q,t,\tau)\,.
\]
It satisfies the following cut-and-join equation
\begin{align*}
 & \frac{\partial F(\mathcal{L};\,q,t,\tau )}{\partial \tau }
=\frac{u}{2}\sum_{\alpha =1}^{L}
    \sum_{i,j\geq 1}
    \bigg(
      ijp^\alpha_{i+j}\frac{\partial^{2}F}
    {\partial p^\alpha_i\partial p^\alpha_j}
             +(i+j)p^\alpha_i p^\alpha_j
        \frac{\partial F}{\partial p^\alpha_{i+j}}
            +ijp^\alpha_{i+j}
      \frac{\partial F}{\partial p^\alpha_i}
        \frac{\partial F}{\partial p^\alpha_j}
    \bigg)\,.
\end{align*}
Restrict the equation to $\vec{\mu}$, we have
\begin{equation}\label{eqn: nonlinear terms}
  \frac{\partial F_{\vec{\mu}}}{\partial \tau }
= \frac{u}{2}
     \bigg(
    \sum_{| \vec{\nu}| =| \vec{\mu}| ,\,\ell(\vec{\nu})=\ell(\vec{\mu})\pm
1}
      \alpha _{\vec{\mu}\vec{\nu}}F_{\vec{\nu}}
    + \textrm{nonlinear terms}
     \bigg) \,,
\end{equation}
where $\alpha _{\vec{\mu}\vec{\nu}}$ is some constant.

Denote by $\deg_u$ the lowest degree of $u$ in a Laurent polynomial of $u$.
The pole structure of $P_{\vec{B}}$ follows from the following degree lemma:

\begin{lemma}\label{lemma: degree}
 $\deg_u F_{\vec{\mu}} \geq \ell(\vec{\mu})-2$.
\end{lemma}

Let $\bigcirc$ be the unknot. Given any
$\vec{A} =(A^1,\ldots,A^L) \in \mathcal{P}^L$, we will obtain the following
limit behavior of quantum group invariants of links at $q\to 1$.
\begin{equation}\label{eqn: normalization of W_A in degree u}
 \lim_{q\rightarrow 1}
 \frac{W_{\vec{A}}(\mathcal{L};\,q,t)}{W_{\vec{A}}(\bigcirc ^{\otimes L};\,q,t)}
=\prod_{\alpha=1}^L
 \xi_{\mathcal{K}_\alpha}(t)^{d_\alpha} \,,
\end{equation}
where $|A^\alpha|=d_\alpha$, $\mathcal{K}_\alpha$ is the
$\alpha$-th component of $\mathcal{L}$,
and $\xi_{\mathcal{K}_\alpha}(t)$, $\alpha=1,\ldots,L$, are
independent of $\vec{A}$.

The following \emph{cut-and-join analysis} will give the desired results.

If $\parallel \vec{\mu}\parallel =1$, it can be verified through the
degree of $u$ in HOMFLY polynomial.

If $\parallel \vec{\mu}\parallel >1$,  by
\eqref{eqn: normalization of W_A in degree u} We can prove that
\[
 \deg_u \Big(\frac{\partial F_{\vec{\mu}}}{\partial \tau} \Big)
=\deg_u F_{\vec{\mu}} \,.
\]
Suppose Lemma \ref{lemma: degree} holds for any
$\parallel \vec{\nu}\parallel <d$.
Now, consider
\[
 S=\{\vec{\nu}:\, |\vec{\nu}| =\vec{d}\}\,,
\]
where $\vec{d}=(d^1,\ldots,d^L)$ and $\sum_{\alpha=1}^L d^\alpha=d$. If
$\vec{\mu}\in S$, one important fact is that the degree of $u$ in the nonlinear
terms in (\ref{eqn: nonlinear terms}) is no less than $\ell(\vec{\mu})-3$. We
will prove Lemma \ref{lemma: degree} holds for $\forall \vec{\mu}\in S$ by
contradiction.

Denote by
\[
 S_{r}
=\{ \vec{\nu}:\vec{\nu}\in S \textrm{ and }\ell(\vec{\nu})=r\}
    \subset S
\]
and
\[
 D_{r}
=\min \{ \deg_{u}F_{\vec{\mu}}:\,\vec{\mu}\in S_r\} \,.
\]

Assume $k$ is the smallest integer such that there exists a
$\vec{\nu}\in S_{k}$ satisfying
\[
\deg _{u}F_{\vec{\nu}}<\ell(\vec{\nu})-2=k-2.
\]
Choose $\vec{\mu}\in S_{k}$ such that $\deg _{u}F_{\vec{\mu}}=D_{k}$.
Rewrite (\ref{eqn: nonlinear terms}) into three parts
\begin{equation}\label{eqn: three parts of log cut-and-join}
 \frac{\partial F_{\vec{\mu}}}{\partial \tau }
=\frac{u}{2}
  \sum_{|\vec{\nu}|=|\vec{\mu}|,\,\ell(\vec{\nu})=\ell(\vec{\mu})-1}
    \alpha _{\vec{\mu}\vec{\nu}} F_{\vec{\nu}}
 +\frac{u}{2}
  \sum_{|\vec{\nu}|=|\vec{\mu}|,\,\ell(\vec{\nu})=\ell(\vec{\mu})+1}
    \alpha _{\vec{\mu}\vec{\nu}} F_{\vec{\nu}}
 +\ast
\end{equation}
where $\ast $ represents the nonlinear part and $\deg _{u}\ast =k-2$.
Here in the equation, $\vec{\nu}$ runs in $S$.

However,
\[
 \deg _{u}\frac{\partial F_{\vec{\mu}}}{\partial \tau}
=\deg_{u}F_{\vec{\mu}}=D_{k}<k-2 \,.
\]
$D_{k-1}\geq (k-1)-2$, so the first part of the sums in the r.h.s of
(\ref{eqn: three parts of log cut-and-join})
is of the lowest order of $u$ at
least $k-2$ and the lowest order of $u$ in
$\frac{\partial F_{\vec{\mu}}}{\partial \tau }$
must come from the summation over $\ell(\vec{\nu})=\ell(\vec{\mu})+1$. Then we
have
\begin{equation}
   1+D_{k+1}\leq D_{k},  \label{ineq: D_k+1 and D_k}
\end{equation}
and
\[
 D_{k+1}\leq D_{k}-1<k-2-1<(k+1)-2\,.
\]
Choose $\vec{\xi}\in S_{k+1}$ such that $\deg _{u}F_{\vec{\xi}}=D_{K+1}$, run
the similar analysis as above. Since (\ref{ineq: D_k+1 and D_k}), one will get
\[
 1+D_{k+2}\leq D_{k+1}.
\]
Thus we obtain
\[
 D_{k}>D_{k+1}>D_{k+2}>...D_{d-1}>D_d \,.
\]

Remember that there is still one last equation we have not used yet. Consider
the cut-and-join equation for the following partition
\[
 \vec{\eta}
=\big( (1^{d_1}),\ldots, (1^{d_L}) \big) \in S_d \,.
\]
Note that this partition has the longest length in $S$, there will be only
the first summand in (\ref{eqn: three parts of log cut-and-join}) left. In
particular, no non-linear terms will appear in the equation. Therefore,
\[
 \frac{\partial F_{\vec{\eta}}}{\partial \tau }
=\frac{u}{2}\sum_{\ell(\vec{\nu})=d-1}
    \alpha _{\vec{\eta}\vec{\nu}}F_{\vec{\nu}} \,.
\]
This implies
\[
 D_d \geq D_{d-1} +1\,,
\]
which contradicts with $D_{d-1}>D_{d}$.This gives a proof of the
existence of \eqref{eqn: P_B}.

Define
        \begin{align*}
            \widetilde{F}_{\vec{\mu}} = \frac{F_{\vec{\mu}} } {\phi_{ \vec{\mu} }(q)}, \qquad \widetilde{Z}_{ \vec{\mu} }
                = \frac{ Z_{ \vec{\mu} } }{ \phi_{ \vec{\mu} }(q) }\,.
        \end{align*}
where
\begin{align*}
    \phi_{\vec{\mu}} (q) = \prod_{\alpha=1}^{L}
        \prod_{j=1}^{\ell(\mu^\alpha)} [\mu^\alpha_j]\,.
\end{align*}
    Lemma \ref{lemma: degree} directly implies the following:
    \begin{corollary} \label{cor: tilde_F_mu}
        For any $\vec{\mu}$, $\widetilde{F}_{\vec{\mu}}$ is of the following form:
        \begin{align*}
            \widetilde{F}_{\vec{\mu}}(q,t)
            = \sum_{\textrm{finitely many } n_\alpha} \frac{ a_\alpha(t) }
                { [n_\alpha]^2 } + \textrm{ polynomial.}
        \end{align*}
    \end{corollary}

Let $\cL$ be the closure of a braid $\beta$ with writhe number $0$.
Cable the $i$-th component of $\beta$, $\beta_i$, by substituting
$\ell(\mu^i)$ parallel strands for each strand of $\beta_i$ and the
 $\ell(\mu^i)$ parallel components are colored by
partitions $(\mu^i_1)$,...,$(\mu^i_{\ell(\mu^i)})$. We take the
closure of this new braid and obtain a new link, denoted by
$\cL_{\vec{\mu}}$.

Let
$\mu^i = (\mu^i_1, \cdots, \mu^i_{\ell_i})$.
We use symbol
\begin{align*}
    \hat{Z}_{ \vec{\mu} } =  Z_{ \vec{\mu} } \cdot \mathfrak{z}_{\vec{\mu}} ; \qquad
        \hat{F}_{\vec{\mu}} = F_{\vec{\mu}} \cdot \mathfrak{z}_{ \vec{\mu} }\,.
\end{align*}
Notice the following fact from the definition of quantum group
invariants:
\begin{align}
    \hat{Z}_{(\mu^1,\cdots, \mu^L)} (\mathcal{L})
    &=\hat{Z}_{(\mu^1_1), \cdots, (\mu^1_{\ell_1}), \cdots, (\mu^L_1),\cdots,
        (\mu^L_{\ell_L})} (\mathcal{L}_{\vec{\mu}}) \,,
        \label{eqn: colored homfly and calbing}
\\
 \hat{F}_{(\mu^1,\cdots, \mu^L)} (\mathcal{L})
    &=\hat{F}_{(\mu^1_1), \cdots, (\mu^1_{\ell_1}), \cdots, (\mu^L_1),\cdots,
        (\mu^L_{\ell_L})} (\mathcal{L}_{\vec{\mu}})
        \label{eqn: connected colored homfly and calbing} \,.
\end{align}

Consider $\delta_n = \sigma_1\cdots \sigma_{n-1}$.
    Let $\fS_A$ be the minimal projection of the Hecke algebra $\mathcal{H}_n \rightarrow \mathcal{H}_A$, and
        \begin{align*}
            \fP_\mu = \sum_A \chi_A(C_\mu) \fS_A\,.
        \end{align*}
        We will apply a lemma of Aiston-Morton \cite{AM} in the following computation:
        \begin{align*}
            \delta_n^n \fS_A = q^{\frac12 \kappa_A} \fS_A\,.
        \end{align*}
Let
        \begin{align*}
            \vec{d} = \big( (d_1),\ldots, (d_L) \big), \qquad \frac{1}{\vec{d}} = \Big( \frac{1}{d_1}, \ldots, \frac{1}{d_L} \Big)\,.
        \end{align*}
        Due to the cabling formula to the length of partition
        \eqref{eqn: colored homfly and calbing} and
        \eqref{eqn: connected colored homfly and calbing},
        we can simply deal with all the color of one row without loss of generality.
        Take framing $\tau_\alpha=n_\alpha + \frac{1}{d_\alpha}$ and choose a braid group representative of $\cL$ such that the writhe number of $\cL_\alpha$ is $n_\alpha$. Denote by $\vec{\tau} = (\tau_1,\ldots,\tau_L)$,
        \begin{align*}
            \hZ_{\vec{d}} \Big(\cL;q,t; \vec{\tau} \Big)
                &= \sum_{\vec{A} } \chi_{\vec{A}} (C_{\vec{d}} )W_{\vec{A}}(\cL;q,t)
                    q^{\frac12 \sum_{\alpha=1}^L
                    \kappa_{A^\alpha} \big(n_\alpha + \frac{1}{ d_\alpha}\big) } \\
                &= t^{\frac12 \sum_\alpha d_\alpha n_\alpha}
                    \Tr \Big( \cL_{\vec{d}}
                    \sum_{\vec{A}} \chi_{\vec{A}}(C_{\vec{d}}) q^{\frac12 \sum_{\alpha} \kappa_\alpha \frac{1}{ d_\alpha} }
                    \bigotimes_\alpha\fS_{A^\alpha} \Big) \\
                &= t^{\frac12\sum_\alpha d_\alpha n_\alpha}
                    \Tr\Big( \cL_{\vec{d}}
                    \sum_{\vec{A}} \chi_{\vec{A}} (C_{\vec{d}}) (\delta_{d_1}\otimes\cdots\otimes \delta_{d_L})                \bigotimes_\alpha\fS_{A^\alpha} \Big)\\
                &= t^{\frac12\sum_\alpha d_\alpha n_\alpha}
                    \Tr \Big( \cL_{\vec{d}} \cdot \otimes_{\alpha=1}^L\delta_{d_\alpha} \cdot
                    \fP_{(1)}^{(d_1)} \otimes \cdots \otimes \fP_{(1)}^{(d_L)} \Big) \,.
        \end{align*}
        Here, $\fP_{(1)}^{(d_\alpha)}$ means that in the projection, we use $q^{d_\alpha}, t^{d_\alpha}$ instead of using $q,t$.
        If we denote by
        \begin{align*}
            \cL\ast Q_{\vec{d}} = \cL_{\vec{d}} \cdot \delta_{\vec{d}} \cdot
                    \fP_{(1)}^{(d_1)} \otimes \cdots \otimes \fP_{(1)}^{(d_L)} \,,
        \end{align*}
        we have
        \begin{align}\label{eqn: Z_mu as homfly of twist}
            \hZ_{\vec{d}} \Big(\cL;q,t; \frac{1}{\vec{d}} \Big) = \mathcal{H}(\cL\ast Q_{\vec{d}})\,.
        \end{align}
        Here $\mathcal{H}$ is the homfly polynomial which is normalized as
        \begin{align*}
            \mathcal{H} (\mathrm{unknot}) = \frac{ t^{\frac12} - t^{-\frac12} }{ q^{\frac12} -q^{-\frac12} }\,.
        \end{align*}
        With the above normalization, for any given link $\cL$, we have
        \begin{align*}
            [1]^L \cdot\mathcal{H} (\cL) \in \mathbb{Q} \big[ [1]^2, t^{\pm\frac12} \big]\,.
        \end{align*}
        Substituting $q$ by $q^{d_\alpha}$ in the corresponding component, it leads to the following:
        \begin{align}\label{eqn: pole of Z_mu}
            \prod_{\alpha=1}^L [d_\alpha] \cdot
            \hat{Z}_{\vec{d}} (\cL; q,t; \vec{\tau})
            \in \mathbb{Q}\big[ [1]^2, t^{\pm\frac12} \big]\,.
        \end{align}

        On the other hand, given any frame $\vec{\omega}=(\omega_1,\ldots,\omega_L)$,
        \begin{align*}
            \hat{Z}_{\vec{\mu}} (\cL;\, q,t; \vec{\omega})
            &= \sum_{\vec{A}} \chi_{\vec{A}} (C_{\vec{\mu}})
                W_{\vec{A}}(\cL;\,q,t;\vec{\omega}) \\
            &= \sum_{\vec{A}} \chi_{\vec{A}}(C_{\vec{\mu}})
                \sum_{\vec{\nu}} \frac{ \chi_{\vec{A}}(C_{\vec{\nu}})} {\mathfrak{z}_{\vec{\nu}}} \hat{Z}_{\vec{\nu}} (\cL;\,q,t)
                q^{\frac12\sum_\alpha \kappa_{A^\alpha}\omega_\alpha }\,.
        \end{align*}
        Exchange the order of summation, we have the following {\em convolution formula}:
        \begin{align}\label{eqn: convolution formula}
            \hat{Z}_{\vec{\mu}}(\cL;\,q,t;\vec{\omega})
            = \sum_{\vec{\nu}} \frac{ \hat{Z}_{\vec{\nu}}(\cL;\,q,t)}
                {\mathfrak{z}_{\vec{\nu}}}
                \sum_{\vec{A}} \chi_{\vec{A}}(C_{\vec{\mu}})
                \chi_{\vec{A}}(C_{\vec{\nu}})
                q^{\frac12\sum_\alpha \kappa_{A^\alpha} \omega_\alpha}\,.
        \end{align}

        This property holds for arbitrary choice of $n_\alpha$, $\alpha=1,\ldots,L$. The coefficients of possible other poles vanish for arbitrary integer $n_\alpha$. $q^{n_\alpha}$ is involved through certain polynomial relation, which implies the coefficients for other possible poles are simply zero. Therefore,
        \eqref{eqn: pole of Z_mu} holds for any frame.

        Now instead of canceling all the poles of $\hat{Z}(\cL)$ according to
        \eqref{eqn: pole of Z_mu}, we consider
        $[c]\hat{Z}_{(c),\vec{d}}$, which is equivalent to consider $[c]^2 \widetilde{Z}_{(c),\vec{d}}$.
        Let $\mathcal{K}$ be the knot labeled by $[c]$. Multiplying $\widetilde{Z}$ by $[c]^2$ cancels all the poles related to $\mathcal{K}$
        according to \eqref{eqn: Z_mu as homfly of twist}. Therefore, $[c]^2\widetilde{Z}_{(c),\vec{d}}$ has the same principle part as $\widetilde{Z}_{\vec{d}}$ except for multiplying by an element in
        $\mathbb{Q}\big[ [1]^2,t^{\pm\frac12}\big]$. Note that
        \begin{align*}
            \widetilde{Z}_{\vec{\mu}} = \sum_{\Lambda\vdash \vec{\mu}}
                \frac{ \widetilde{F}_{\Lambda} }{ \Aut |\Lambda| }\,.
        \end{align*}
        $[c]^2(\widetilde{Z}_{(c),\vec{d}}-\widetilde{F}_{(c),\vec{d}})$ contains all the principle terms from $\cL$ omitting $\mathcal{K}$. Therefore,
        \begin{align*}
            [c]^2 \widetilde{F}_{(c),\vec{d}} \in \mathbb{Q}\big[ [1]^2, t^{\pm\frac12} \big]\,.
        \end{align*}
        In the above discussion, $\mathcal{K}$ can be chosen to be any component of $\cL$. We thus proved the following proposition:

        \begin{proposition} \label{prop: Z_mu pole structure}
            Notation as above, we have:
        \begin{align}
            \prod_{\alpha=1}^L [d_\alpha] \cdot
                \hZ_{\vec{d}} (\cL;q,t)
            & \in \mathbb{Q} [ [1]^2, t^{\pm\frac12}]; \\
                [d_\alpha]^2 \widetilde{F}_{\vec{d}}
                ( \cL; q,t )
            &\in \mathbb{Q}\big[ [1]^2, t^{\pm\frac12} \big], \forall \alpha\,.
        \end{align}
        \end{proposition}

        Similarly, We can obtain
        \begin{align} \label{eqn: pole of tilde_F_mu}
            \widetilde{F}_{\vec{d}}(\cL;q,t)
            = \frac{ H_{\vec{d}/D_{\vec{d}}} ( t^{D_{\vec{d}}} ) } {D_{\vec{d}}\cdot [D_{\vec{d}}]^2} +
            \textrm{polynomial in $[D_{\vec{d}}]^2$ and $t^{\pm\frac12 D_{\vec{d}}}$}\,.
        \end{align}
        Once again, due to arbitrary choice of $n_\alpha$, we know the above pole structure of $\widetilde{F}_{\vec{d}}$ holds
        for any frame.

\begin{proposition} \label{prop: pole structure of tilde_F_mu}
            Notations are as above. Assume $\cL$ is labeled by the color $\vec{\mu}=(\mu^1,\ldots,\mu^L)$. Denote by $D_{\vec{\mu}}$ is the greatest
            common divisor of $\{\mu^1_1$, $\ldots$ ,$\mu^1_{\ell(\mu^1)}$,$\ldots$, $\mu^i_j,\ldots, \mu^L_{\ell(\mu^L)}\}$.
            $\widetilde{F}_{\vec{\mu}}$ has the following structure:
            \begin{align*}
                \widetilde{F}_{ \vec{\mu} } (q,t)  = \frac{ H_{\vec{\mu}/D_{\vec{\mu}}}( t^{D_{\vec{\mu}}} ) } 
                	{D_{\vec{\mu}}\cdot [D_{\vec{\mu}} ]^2 } +
                f(q,t)\,,
            \end{align*}
            where $H_{\vec{\mu}/D_{\vec{\mu}}}(t)$ only depends on $\vec{\mu}/D_{\vec{\mu}}$ and $\cL$, $f(q,t)\in\mathbb{Q}\big[ [1]^2, t^{\pm\frac12} \big]$.
        \end{proposition}

        \begin{remark}
            In Proposition \ref{prop: pole structure of tilde_F_mu}, it is very interesting to interpret in topological string side that
            $H_{\vec{\mu}/D_{\vec{\mu}}}(t)$ only depends on $\vec{\mu}/D_{\vec{\mu}}$ and $\cL$. The principle term is generated due to
            summation of counting rational curves and independent choice of $k$ in the labeling color
            $k\cdot\vec{\mu}/D_{\vec{\mu}}$. This phenomenon simply tells us that contributions of counting rational curves in
            the labeling color $k\cdot\vec{\mu}/D_{\vec{\mu}}$ are through multiple cover contributions of
            $\vec{\mu}/D_{\vec{\mu}}$.
        \end{remark}

\subsection{Integrality}

By the definition of $P_{\vec{B}}(q,t)$, comparing with Proposition
\ref{prop: pole structure of tilde_F_mu}, we have the following
computation:
\begin{align*}
        P_{\vec{B}}(q,t)
            &=\sum_{\vec{\mu}} \frac{\chi_{\vec{B}}(\vec{\mu})}{\phi_{\vec{\mu}}(q)}
                \sum_{d|\vec{\mu}} \frac{\pi (d)}{d}F_{\vec{\mu}/d}(q^{d},t^{d}) \\
            &=\sum_{ \vec{\mu} } \chi_{\vec{B}} ( \vec{\mu} ) \sum_{d| \vec{\mu} }
                \frac{ \mu(d) }{d} \widetilde{F}_{\vec{\mu}/d} (q^d, t^d) \\
            &=\sum_{ \vec{\mu} } \chi_{\vec{B}} ( \vec{\mu} ) \sum_{d| D_{\vec{\mu}} }\frac{ \mu(d) }{d}
                \frac{ H_{\vec{\mu}/D_{\vec{\mu}} }(t^{D_{\vec{\mu}}}) }
                {D_{\vec{\mu}/d}\cdot [D_{\vec{\mu}}]^2 } + \textrm{polynomial} \\
            &=\sum_{ \vec{\mu} } \chi_{\vec{B}} ( \vec{\mu} ) \delta_{1, D_{\vec{\mu}}}
                \frac{ H_{\vec{\mu}/D_{\vec{\mu}} }(t^{D_{\vec{\mu}}}) }
                {D_{\vec{\mu}}\cdot [D_{\vec{\mu}}]^2 } + \textrm{polynomial}\,,
    \end{align*}
    where $\delta_{1,n}$ equals $1$ if $n=1$ and $0$ otherwise.
    It implies that $P_{\vec{B}}$ is a rational function which only has pole at $q=1$.
    In the above computation, we used a fact of M\"obius inversion,
    \begin{align*}
        \sum_{d|n} \mu(d) = \delta_{1,n}\,.
    \end{align*}

Therefore, for each $\vec{B}$,
\[
 \sum_{g= 0}^\infty \sum_{Q\in \mathbb{Z}/2}
    N_{\vec{B};\,g,Q}(q^{-1/2}-q^{1/2})^{2g}t^Q
\in \mathbb{Q}[(q^{-1/2}-q^{1/2})^2, t^{\pm \frac{1}{2}}]\,.
\]

The integrality of $N_{\vec{B};\,g,Q}$ can be derived from the
following theorem:

\begin{Theorem}\label{thm: integrality}
 We have
\begin{align}
 \sum_{g=0}^\infty \sum_{Q\in\mathbb{Z}} N_{\vec{B};\,g,Q}
    (q^{-1/2}-q^{1/2})^{2g} t^Q
\in \mathbb{Z}[(q^{-1/2}-q^{1/2})^2, t^{\pm 1/2}]\,.
\end{align}
\end{Theorem}

\begin{remark}
 The above theorem also implies the refined integral invariants
$N_{\vec{B};\,g,Q}$ vanish at large genera (also for large $Q$).
\end{remark}

Define $y=(y^1,y^2,\ldots,y^L)$ where
\[
 y^\alpha = (y^\alpha_1,y^\alpha_2,\ldots)
\]
is a set of arbitrarily chosen variables. Denote by $\Omega(y)$ the set of
all symmetric function in $(y^1,y^2,\ldots,y^L)$ with integral coefficients.

We start by defining the following special ring which characterizing the
algebraic structure of Chern-Simons partition function.
\begin{align*}
  \mathfrak{R}(y;q,t)
= \bigg\{ \frac{a(y;q,t)}{b(q)}:\, a(y;q,t)\in
    \Omega(y)[[1]_q^2,t^{\pm 1/2}],\,
    b(q)=\prod_{n_k} [n_k]^2_q\in \mathbb{Z}[[1]_q^2]
   \bigg\}  \,.
\end{align*}
Given $\frac{f(y;\,q,t)}{b(q)}\in \mathfrak{R}(y;\,q,t)$, if $f(y;\,q,t)$ is a
primitive polynomial in terms of $q^{\pm 1/2}$, $t^{\pm 1/2}$ and Schur
functions of $y$, we call $\frac{f(y;\,q,t)}{b(q)}$ is primitive.

Given any $\frac{r}{s}h(y;\,q,t)$ where $h(y;\,q,t)\in\mathfrak{R}(y;\,q,t)$ is
primitive, define
\[
 \Ord_p\Big(\frac{r}{s}h(y;\,q,t)\Big)
= \Ord_p \Big( \frac{r}{s}\Big)
\]
for any prime number $p$.

We will consider the following generating series
\[
 T_{\vec{d}} = \sum_{|\vec{B}|=\vec{d}} s_{\vec{B}}(y) P_{\vec{B}}(q,t)
\]
After some calculations, we have
\begin{equation}\label{eqn: integrality eqn}
 T_{\vec{d}}
=
 q^{| \vec{d}| }
    \sum_{k|\vec{d}}\frac{\mu (k)}{k}
    \sum_{\mathfrak{A}\in \mathcal{P(P}^{n}),
        \parallel\mathfrak{A}\parallel =\vec{d}/k}
    \theta_{\mathfrak{A}} \prod_{j=1}^{\ell(\mathfrak{A})}
    W_{\mathfrak{A}_j}(q^{k},t^{k})
    s_{\mathfrak{A}_j}((z)^{k})
\end{equation}

\begin{proposition}
\label{prop: T_d in L}$T_{\vec{d}}(y;\,q,t)\in \mathfrak{R}(y;v,t)$.
\end{proposition}

Proposition \ref{prop: T_d in L} implies $q=0$ is a pole of $T_{\vec{d}}$.
Since $T_{\vec{d}}$ can be written as
\[
 T_{\vec{d}}
=\sum_{g= 0}^\infty \sum_{Q\in \mathbb{Z}/2}
    \bigg( \sum_{|\vec{B}|=\vec{d}}N_{\vec{B};\,g,Q}s_{\vec{B}}(y)
    \bigg)
    (q^{-1/2}-q^{1/2})^{2g-2}t^{Q}
\]
by the existence of \eqref{eqn: P_B}.
If there are infinitely many $N_{\vec{B};\,g,Q}$ nonzero, $q=0$ is then an
essential singularity point of $T_{\vec{d}}$, which is a contradiction.
Therefore, for each $\vec{B}$,
\[
 \sum_{g= 0}^\infty \sum_{Q\in \mathbb{Z}/2}
    N_{\vec{B};\,g,Q}(q^{-1/2}-q^{1/2})^{2g}t^Q
\in \mathbb{Q}[(q^{-1/2}-q^{1/2})^2, t^{\pm \frac{1}{2}}]\,.
\]

On the other hand, by Proposition \ref{prop: T_d in L},
$T_{\vec{d}}\in \mathfrak{R}(y;\,v,t)$ and
$\Ord_p T_{\vec{d}} \geq0$ for any
prime number $p$. We have
\[
 \sum_{|\vec{B}|=\vec{d}} N_{\vec{B};\,g,Q} s_{\vec{B}} (y)
\in \Omega(y) \,,
\]
which implies $N_{\vec{B};\,g,Q}\in\mathbb{Z}$.

Combining the above discussions, we have
\[
    \sum_{g= 0}^\infty \sum_{Q\in \mathbb{Z}/2}
    N_{\vec{B};\,g,Q}(q^{-1/2}-q^{1/2})^{2g}t^{Q}
\in \mathbb{Z}[(q^{-1/2}-q^{1/2})^{2},t^{\pm 1/2}] \, .
\]

To prove Proposition \ref{prop: T_d in L}, we combine the study of
multi-cover contribution and $p$-adic argument. For any give prime
number $p$, the following observation is important for the $p$-adic
argument
\[
 \bigg\{ \mathfrak{B}:\; \sum_{j=1}^{\ell(\mathfrak{B})} |\mathfrak{B}_j|
    = p\,\vec{d}
    \textrm{ and } \Ord_p(\theta_{\mathfrak{B}})<0
 \bigg\}
=\bigg\{ \mathfrak{A}^{(p)}:\; \sum_{j=1}^{\ell(\mathfrak{A})}
    |\mathfrak{A}_j| = \vec{d}\,
 \bigg\}\,.
\]
Matching the following terms, finally we have
\begin{equation}\label{eqn: calculation for lemma of multicover}
  \Ord_{p}\Big( \theta _{\mathfrak{A}^{(p)}}W_{\mathfrak{A}^{(p)}}(q,t)s_{%
    \mathfrak{A}^{(p)}}(z)-\frac{1}{p}\theta _{\mathfrak{A}}W_{\mathfrak{A}%
    }(q^{p},t^{p})s_{\mathfrak{A}}(z^{p})\Big)
\geq 0 \,.
\end{equation}
Let
\[
 \Phi_{\vec{d}}\left( y;\,q,t\right)
=\sum_{\mathfrak{A}\in \mathcal{P(P}^{n}),\,
 \parallel\mathfrak{A}\parallel =\vec{d}}\theta _{\mathfrak{A}}
 W_{\mathfrak{A}}(q,t)s_{\mathfrak{A}}(z) \, .
\]
The following inequality can be obtained from
\eqref{eqn: calculation for lemma of multicover}:
\begin{align*}
 \Ord_p \big( \Phi_{p\vec{d}}\,(y;\,q,t )
 -\frac{1}{p} \Phi_{\vec{d}}\, ( y^{p};q^{p},t^{p}) \big)
\geq 0 \,.
\end{align*}
Therefore,
\begin{align*}
    T_{\vec{d}}
&= q^{|\vec{d}|} \sum_{k|\vec{d}} \frac{\mu(k)}{k}
    \sum_{
      \mathfrak{A}\in\mathcal{P(P}^L),\,
      \parallel\mathfrak{A}\parallel =\vec{d}/k
     }
    \theta_{\mathfrak{A}} W_{\mathfrak{A}}(q^{k},t^{k})
    s_{\mathfrak{A}}(z^{k})
\\
&= q^{|\vec{d}|}
    \sum_{k|\vec{d},\, p\nmid k} \frac{\mu(k)}{k}
    \Big(
        \Phi_{\vec{d}/k}(y^{k};\,q^{k},t^{k})
    -\frac{1}{p}\Phi _{\vec{d}/(pk)}(y^{pk};\,q^{pk},t^{pk})
    \Big) \,.
\end{align*}
This implies
\begin{equation}\label{ineq: Ord_p T_d geq 0}
 \Ord_{p}T_{\vec{d}}\geq 0\,.
\end{equation}
Combining with \eqref{eqn: P_B}, we thus prove Proposition
\ref{prop: T_d in L}.

\section{Concluding remarks}

In this section, we briefly discuss some interesting problems
related to string duality which may be approached through the
techniques developed in this paper.

 Let
$
 \mathbf{p}
=(\mathbf{p}^1, \ldots, \mathbf{p}^L), $ where
$
 \mathbf{p}^\alpha
=(p^\alpha_1,p^\alpha_2, \ldots,) \,. $ Defined the following
generating series of open Gromov-Witten invariants
\[
 F_{g,\vec{\mu}}(t,\tau)
=\sum_{\beta} K^\beta_{g,\vec{\mu}}(\tau) e^{\int_\beta \omega}
\]
where $\omega$ is the K\"ahler class of the resolved conifold,
$\tau$ is the framing parameter and
\begin{align*}
 t=e^{\int_{\mathbb{P}^1} \omega}, \quad \textrm{and} \quad
    e^{\int_\beta \omega} = t^Q\,.
\end{align*}
Consider the following generating function
\[
 F(\mathbf{p};\,u,t;\,\tau)
=\sum_{g=0}^\infty \sum_{\vec{\mu}} u^{2g-2+\ell(\vec{\mu})}
    F_{g,\vec{\mu}}(t;\, \tau)
    \prod_{\alpha=1}^L p^\alpha_{\mu^\alpha} \,.
\]
It satisfies the log cut-and-join equation
\[
 \frac{\partial F(\mathbf{p};\,u,t;\,\tau)}{\partial \tau}
=\frac{u}{2} \sum_{\alpha=1}^L \mathfrak{L}_\alpha
    F(\mathbf{p};\,u,t;\,\tau)\,.
\]
Therefore, duality between Chern-Simons theory and open
Gromov-Witten theory reduces to verifying the uniqueness of the
solution of cut-and-join equation.

Cut-and-join equation for Gromov-Witten side comes from the
degeneracy and gluing procedure while uniqueness of cut-and-join
system should in principle be obtained from the verification at some
initial value. However, it seems very difficult to find a suitable
initial value. A new hope might be found in our development of
cut-and-join analysis. In the log cut-and-join equation the
non-linear terms reveals the important recursion structure. For the
uniqueness of cut-and-join equation, it will appear as the vanishing
of all non-linear terms. We will put this in our future research.

Our proof of the LMOV conjecture has shed new light on the famous
volume conjecture. See for example the discussions in \cite{Gukov}.
The cut-and-join analysis we developed in this paper combined with
rank-level duality in Chern-Simons theory seems to provide a new way
to prove the existence of the limits of quantum invariants.

There are also other open problems related to LMOV conjecture. For
example, quantum group invariants satisfy skein relation which must
have some implications on topological string side as mentioned in
\cite{LM}. One could also rephrase a lot of unanswered questions in
knot theory in terms of open Gromov-Witten theory. We hope that the
relation between knot theory and open Gromov-Witten theory will be
explored much more in detail in the future. This will definitely
open many new avenues for future research.

\bibliographystyle{mrl}

\end{document}